
\magnification1200
\input amstex.tex
\documentstyle{amsppt}
\nopagenumbers
\hsize=12.5cm
\vsize=18cm
\hoffset=1cm
\voffset=2cm

\footline={\hss{\vbox to 2cm{\vfil\hbox{\rm\folio}}}\hss}

\def\DJ{\leavevmode\setbox0=\hbox{D}\kern0pt\rlap
{\kern.04em\raise.188\ht0\hbox{-}}D}

\def\txt#1{{\textstyle{#1}}}
\baselineskip=13pt
\def\hf{{\textstyle{1\over2}}}
\def\b{\beta}
\def\d{{\,\roman d}}
\def\e{\varepsilon}
\def\E{{\roman e}}

\def\b{\beta} \def\g{\gamma}
\def\G{\Gamma}

\def\t{\theta}
\def\={\;=\;}

\def\zt{\zeta(\hf+it)}

\def\r{\rho}
 
\def\z{\zeta}

\def\t{\theta}
\def\hf{{\textstyle{1\over2}}}
\def\txt#1{{\textstyle{#1}}}

\def\le{\leqslant} \def\ge{\geqslant}
\font\tenmsb=msbm10
\font\sevenmsb=msbm7
\font\fivemsb=msbm5
\newfam\msbfam
\textfont\msbfam=\tenmsb
\scriptfont\msbfam=\sevenmsb
\scriptscriptfont\msbfam=\fivemsb
\def\Bbb#1{{\fam\msbfam #1}}

\def \NN {\Bbb N}
\def \CC {\Bbb C}
\def \RR {\Bbb R}
\def \ZZ {\Bbb Z}

\font\ff=cmr8
\baselineskip=13pt

\font\teneufm=eufm10
\font\seveneufm=eufm7
\font\fiveeufm=eufm5
\newfam\eufmfam
\textfont\eufmfam=\teneufm
\scriptfont\eufmfam=\seveneufm
\scriptscriptfont\eufmfam=\fiveeufm
\def\mathfrak#1{{\fam\eufmfam\relax#1}}

\font\tenmsb=msbm10
\font\sevenmsb=msbm7
\font\fivemsb=msbm5
\newfam\msbfam
     \textfont\msbfam=\tenmsb
      \scriptfont\msbfam=\sevenmsb
      \scriptscriptfont\msbfam=\fivemsb
\def\Bbb#1{{\fam\msbfam #1}}

\def \NN {\Bbb N}
\def \CC {\Bbb C}

\def \RR {\Bbb R}
\def \ZZ {\Bbb Z}

  \def\rightheadline{{\hfil{\ff
  A note on the zeros of the zeta-function}\hfil\tenrm\folio}}

  \def\leftheadline{{\tenrm\folio\hfil{\ff
   Aleksandar Ivi\'c }\hfil}}
  \def\emptyheadline{\hfil}
  \headline{\ifnum\pageno=1 \emptyheadline\else
  \ifodd\pageno \rightheadline \else \leftheadline\fi\fi}

\topmatter
\title
A note on the zeros of the zeta-function of Riemann
\endtitle
\author   Aleksandar Ivi\'c
 \endauthor

\nopagenumbers

\medskip

\address
Aleksandar Ivi\'c, Serbian Academy of Science and Arts, Knez Mihailova 35,
11000 Beograd, Serbia \medskip
\endaddress
\keywords Riemann zeta-function, Hardy's function $Z(t)$, definition of $N(T)$,
Riemann hypothesis
\endkeywords
\subjclass
11M06  \endsubjclass

\bigskip
\email {
\tt
aleksandar.ivic\@rgf.bg.ac.rs, aivic\_2000\@yahoo.com }\endemail
\dedicatory
\enddedicatory
\abstract
It is shown explicitly how the sign of Hardy's function $Z(t)$ depends on the parity
of the zero-counting function $N(T)$. Two existing definitions of this function are
analyzed, and some related problems are discussed.

\endabstract
\endtopmatter

\document

\medskip
\head
1. Introduction
\endhead

\medskip
Let the Riemann zeta-function be, as usual,
$$
\z(s) \;=\; \sum_{n=1}^\infty n^{-s}\qquad(\Re s > 1).
$$
For $\Re s \le 1$ one defines $\z(s)$ by analytic continuation (see the monographs
of H.M. Edwards [1], the author  [2], [3] and E.C. Titchmarsh
[9] for an extensive account on $\z(s)$).
\medskip
Hardy's function (see the author's
monograph [5] for its properties) is
$$
Z(t)  := \zt\bigl(\chi(\hf+it)\bigr)^{-1/2}, \;\z(s) \;=\; \chi(s)\z(1-s),\leqno(1.1)
$$
so that
$$
\chi(s) = \frac{\G(\hf(1-s))}{\G(\hf s)}\pi^{s-1/2}.\leqno(1.2)
$$
The relation $\z(s) = \chi(s)\z(1-s)$ holds $\forall s \in \CC$. It was proved first by B. Riemann in 1859,
and this is the well-known functional equation for the zeta-function. In 1914 G.H. Hardy introduced $Z(t)$
to prove that there are infinitely many zeros of $\z(s)$  on the so-called ``critical line''
$\Re s = 1/2$.

\medskip
From (1.1), (1.2)  and $\chi(s)\chi(1-s) = 1$ follows that $Z(t)$ is a smooth, even,
real-valued function of the real variable
$t$, for which $|Z(t)|=|\zt|$. Further we have
$$
{(\chi(\hf+it)\bigr)}^{-1/2} =  \pi^{-it/2}{\G^{1/2}({1\over4}+\hf it)\over
\G^{1/2}({1\over4}-\hf it)} = \pi^{-it/2}\frac
{\G({1\over4}+\hf it)}{|\G({1\over4}+\hf it)|}
 := {\roman e}^{i\t(t)},\leqno(1.3)
$$
say. Thus (1.3) implies that, for $t \in \RR$,
$$
\eqalign{
\t(t) &=  -{1\over2i}\log\chi(\hf+it) = \frac{1}{i}\left\{\log \G(\txt{1\over4} + \hf it) -
\log |\G(\txt{1\over4} + \hf it)|\right\} - \hf t\log\pi\cr&
= \arg \G(\txt{1\over4} + \hf it) - \hf t\log\pi =
\Im\left\{\log\Gamma({\txt{1\over 4}} + \hf it)\right\}
- \hf t\log\pi\;\in \RR.\cr}
\leqno(1.4)
$$

\medskip
The  problem of the distribution of complex zeta zeros $\rho_n = \b_n + i\g_n$ is a fundamental
one in zeta-function theory. In view of the functional equation and $\overline{ \z(s)} = \zeta({\bar s})$,
we may assume without loss of generality that $\hf \le \b_n< 1$ (since it is elementary that there
are no zeros in the region $\Re s \ge1$) and $\g_n > 0$. Here  the ordinates of zeros are ordered
as $0 < \g_1 \le \g_2 \le \g_3 \le \cdots$, with ordinates belonging to multiple zeros (if any) being
considered as different. All known zeros are simple, and the first three in the upper
complex half-plane  are approximately
$$
 {1\over2} + i14.134725\ldots\,,\;{1\over2} + i21.022039\ldots\,,\;
{1\over2} + i25.010857\ldots\,.
$$
For recent results on zeta-zeros, see e.g., S. Wedeniwski [10].

\medskip
If $\g\;(>0)$ denotes generic ordinates of zeta zeros, then  the integer-valued function
$$
N(T) \;:=\; \sum_{0<\g\le T}1\leqno(1.5)
$$
counts the number of $\g$'s not exceeding $T$ (with multiplicities included),
and is therefore well-defined
for any $T>0$. This is the definition from E.C. Titchmarsh [9], p. 211. There is another definition
of $N(T)$, namely
$$
\eqalign{
N(T) &\;:=\; \sum_{0<\g< T}1\qquad(T\ne \g),\cr
N(T) &\;:=\;\hf\bigl(N(T-0) + N(T+0)\bigr)\qquad(T = \g).\cr}\leqno(1.6)
$$
This definition is already in the title  of A. Selberg's paper  [8], but I am not certain in which works the above
definitions appeared first. Unless $T = \g$ both definitions of $N(T)$ are equal, and in the applications which
I have seen so far it does not matter which one of them is used. However, if $N(\g-0)$ and $N(\g+0)$
are of different parity, then according to Selberg's definition $N(\g)$ is
of the form $2k - \hf, k\in\NN$, and
thus not an integer. We shall see later that this fact is of some significance.

\medskip

A version of the well known Riemann - von Mangoldt formula for $N(T)$
(see p. 212 of E.C. Titchmarsh [9] and equation (1.45) of [5]) asserts that
$$
N(T) = \frac{1}{\pi}\,\t(T) + 1 + S(T)\qquad(T \ne \g), \leqno(1.7)
$$
with
$$
S(T) \;:=\; \frac{1}{\pi}\arg \z(\hf + iT)\qquad(T \ne \g).\leqno(1.8)
$$
Here $\t(T)$ is as in (1.4), while $\arg\z(\hf+iT)$ is evaluated by continuous variation starting from
$\arg \z(2)=0$ and
proceeding along straight lines, first up to $2+iT$ and then to $1/2+iT$,
assuming that $T$ is not an ordinate of a zeta zero. If $T$ is an ordinate
of a zero, then we set $S(T) = S(T+0)$. Again, this is according to [9], while Selberg [8] defines
$$
S(T) \;=\; \hf\bigl(S(T-0) + S(T+0)\bigr)\qquad(T=\g),
$$
which is in tune with (1.6). From (1.7) and (1.8) it is not difficult to deduce that
$$
N(T) = {T\over2\pi}\log\Bigl({T\over2\pi}\Bigr) - {T\over2\pi} + {7\over8}
+S(T) + O\Bigl({1\over T}\Bigr),
$$
where the $O$-term is a continuous function of $T$, and $S(T) = O(\log T)$.
\medskip
\head
2. Connection between $Z(t)$ and $N(t)$
\endhead

\medskip
From  $|Z(t)| = |\zt|$ and the fact that $Z(t)$ is real
it follows that either $Z(t) = |\zt|$ or $Z(t) = -|\zt|$, but it is not clear
which of these relations holds.
A new formula s presented here, which settles this question by connecting the 
functions $Z(t)$ and $N(T)$. This is

\medskip
THEOREM 1. {\it If $\,t\,$ is not the ordinate of any zeta-zero $\b_n+i\g_n$, then}
$$
Z(t) = (-1)^{N(t)+1}|\zt|\qquad(t>0).\leqno(2.1)
$$

\medskip
Theorem 1 provides an unconditional connection between the sign of $Z(t)$ and the parity
of the zero-counting function $N(t)$. We assume that $N(T)$ is defined by (1.5).
We start from (1.7), which gives
$$
\t(T) = \pi N(T) - \pi S(T) - \pi\qquad(T\ne \g_n),\leqno(2.2)
$$
where $T$ is not an ordinate of any zeta zero $\b_n+i\g_n$.
We record the explicit representation of $\t(t)$ (see (1.21) and (1.22) of [5]),
which follows from (1.4) and Stirling's classical formula for the gamma-function. This is
$$
\theta(t) = {t\over 2}\log {t\over 2\pi} - {t\over 2} -
{\pi\over 8} + \Delta(t)\qquad(t>0).
$$
Here
$$
\Delta(t) := {t\over 4}\log\left(1 + {1\over 4t^2}\right) +
{1\over 4}\arctan{1\over 2t} + {t\over 2}\int\limits_0^\infty
{\psi(u)\over {(u + {1\over 4})}^2 + {(\hf t)}^2}\d u
$$
with
$$
\psi(x) := x - [x] - \hf = - \sum_{n=1}^\infty{\sin(2n\pi x)\over n\pi}
\qquad(x\not\in\ZZ).
$$
We have the approximation (see e.g., p. 120 of H.M. Edwards [1])
$$
\t(t)  = {t\over2}\log {t\over2\pi} - {t\over2} - {\pi\over8} +
{1\over48t} + {7\over5760t^3}+O\Bigl({1\over t^{5}}\Bigr).
$$

\medskip
For any complex number $z$ one has $z = |z|\E^{i\arg z}$.
Hence from the defining relation (1.1) and (1.3) it follows that
$$
Z(t) = \zt\E^{i\t(t)} = |\zt|\E^{i\arg \z(\frac{1}{2}+it) + i\t(t)}.\leqno(2.3)
$$
We have, on using (1.8) and (2.2),
$$
i\arg \zt + i\t(t) = i\pi S(t) + i\pi N(t) - \pi i S(t) - i\pi = i\pi N(t) - i\pi.
$$
But since $\E^{\pm\pi i} = -1$, we obtain from (2.2)
$$
Z(t) = (-1)^{N(t)+1}|\zt|\qquad(t \ne \g),
$$
as asserted in Theorem 1. Clearly there is no restriction on $t$, thus (2.1) holds for any $t >0$
such that $t\ne \g$, i.e., if $t$ is not an ordinate of any zeta zero. In the course of the proof
we used  (1.5). However, in the case when $t\ne \g$, which
is assumed in the formulation of Theorem 1, both definitions of $N(T)$ in (1.5) and (1.6)
coincide. Therefore Theorem 1
remains (unconditionally) true regardless of which definition of $N(T)$ one uses.

\medskip
\head
3. Discussion of the two definitions of $N(T)$
\endhead
If the famous Riemann Hypothesis (RH, all complex zeros of $\z(s)$ satisfy $\Re s = 1/2$) is true,
then the (real) zeros of $Z(t)$ correspond to the zeros $1/2+it$ of $\z(s)$. In this case both
sides of (2.1) vanish, so Theorem 1 holds in this case for all $t>0$, regardless of which definition
of $N(T)$ is used.

\medskip
Suppose now that the RH fails. Then there exist $\b,\g\in\RR$ such that $1/2 < \b < 1, \g >0$ and
$\z(\b + i\g) = \z(1-\b+i\g) = 0$. Moreover, the functional equation $\z(s) = \chi(s)\z(1-s)$
implies that both $\b+i\g$ and $1-\b+i\g$ have the same multiplicity. Recall that if $r = m(\rho)$
denotes the multiplicity of the complex zero $\rho = \b + i\g$
of  $\z(s)$, then
$$
\z(\r) = \z'(\r)
= \ldots = \z^{(r-1)}(\r) = 0, \;\roman{but}\; \z^{(r)}(\r) \not = 0.
$$
Furthermore a zero $\rho$ is simple if $m(\rho)=1$, namely if $\z(\rho) = 0$ but $\z'(\rho) \ne0$.
The problem of evaluating the multiplicities of zeta zeros is a deep one (see  the author's
papers [4], [6], and A.A. Karatsuba [7]). Suppose we use the  definition of $N(T)$ given by (1.5).
Note that by continuity $Z(\g) = Z(\g-0)$, and (2.1) holds for $t = \g-\e$, and
sufficiently small $\e>0$.
The contribution of $N(\g-0)$  differs from that of $N(\g)$ by
the contribution of $S := \sum_{j=1}^k 2m(\b_j +i\g)$, where $\hf < \b_1 < \cdots< \b_k<1$ and
$\z(\b_j+i\g) = 0\;(j = 1, \ldots, k)$. Here we assumed that the $\b_j$'s are the real parts of
all different values of $\hf < \b < 1$ such that $\z(\b+i\g)=0$. Whatever the value of $k$ is,
the sum $S$ is an even natural number $2K, \; K = K(k)$, hence $N(\g) = N(\g-0) + 2K$.
Then we have
$$
\eqalign{
Z(\g) &= Z(\g-0) = (-1)^{N(\g-0)+1}|\z(\hf + i(\g-0))| \cr&= (-1)^{N(\g)-2K+1}|\z(\hf+i\g)|
= (-1)^{N(\g)+1}|\z(\hf+i\g)|.\cr}
$$
This shows that (2.1) holds if $t = \g$, and the same conclusion holds if we use $Z(\g) = Z(\g+0)$.

\medskip
However, if we use the second  definition of $N(T)$ in (1.6), then (2.1) does not have to hold
for $t = \g$, no matter how we define $Z(\g)$ by a limit process. Namely, as already noted,
if $N(\g-0)$ and $N(\g+0)$
are of different parity, then according to Selberg's definition $N(\g)$ is
of the form $2k - \hf, k\in\NN$, hence $(-1)^{N(\g)+1}$ is neither +1 nor -1. This happens e.g.,
if $N(\g-0)$ is even, but $\hf + i\g$ is a simple zero, and there are no other zeta zeros $\rho$
with $\Im \rho = \g$. We see that,
in the context of Theorem 1, it does matter which definition of $N(T)$ is used.
If we use the  definition of $N(T)$ in (1.5), then actually (2.1) holds for all $t>0$ unconditionally.

\medskip
In fact, the argument of the previous section provides another proof that (2.1) holds unconditionally
if (1.5) is used.  In view of (1.2) one has
$$Z(0) = \chi^{-1/2}(\hf)\zeta(\hf) = \zeta(\hf) < 0,
$$
which establishes (2.1) for $0 < t \le \g_1$.
We can then easily verify that (2.1) holds for $t < 10^{13}$, e.g.,
when $t$ does not exceed the ordinate of
the largest known zero (see [10]). Since all known zeros are simple and satisfy the RH, $N(T)$ jumps by
unity at each $\g$, so (2.1) holds. Let $\g_N$ be the largest known ordinate of a zeta-zero, so that
(2.1) holds for $t \le \g_N$. Since $\hf + i\g_N$ is a simple zero, then $Z(t)$ changes sign at $\g_N$,  
so (2.1) holds also in $[\g_N, \g_{N+1})$. At $t=\g_{N+1}$  two cases are possible.
The first case is when $\hf + i\g_{N+1}$ is a zeta-zero, and there is no $\b > \hf$
such that $\z(\b + i\g_{N+1}) =0$. In that case $N(t)$ jumps at $t = \g_{N+1}$
by $m(\hf + i\g_{N+1})$. If this number is odd, then $Z(t)$ changes sign at $\g_{N+1}$ and (2.1)
is true for $\g_{N+1} \le t < \g_{N+2}$. If this number is even, then $Z(t)$ maintains its sign
in the same interval, and again (2.1) is true for $\g_{N+1} \le t < \g_{N+2}$. The same is true
if there are zeros $\rho$ off the critical with this $\g_{N+1}$. The contribution of these zeros
to $N(\g_{N+1})$, as seen in the previous section, is an even number. In the second case,
when $\z(\hf + i\g_{N+1}) \ne0$, but there is a $\b\; (> \hf)$ such that $\z(\b + i\g_{N+1}) =0$
is similar. The contribution to $N(\g_{N+1})$ is an even number, and so
it does not affect the sign of $(-1)^{N(t)+1}$. In any case (2.1) is verified for
$\g_{N+1} \le t < \g_{N+2}$, and the proof of (2.1) may be given inductively.

\medskip
{\bf Corollary 1}. Suppose $Z(t_1) \ne 0, Z(t_2) \ne 0$. Then $Z(t_1)$ and $Z(t_2)$ are both of the
same sign if and only if $N(t_1)$ and $N(t_2)$ are of the same parity, provided that $N(T)$ is
defined by (1.5).

\medskip
This follows from (2.1), since we have
$$
\eqalign{
Z(t_1)Z(t_2) &= (-1)^{N(t_1)+1}(-1)^{N(t_2)+1}|\z(\hf + it_1)||\z(\hf+it_2)|\cr&
= (-1)^{N(t_1)+N(t_2)}|\z(\hf+it_1)\z(\hf+it_2)|.\cr}
\leqno(2.5)
$$
Thus $Z(t_1)$ and $Z(t_2)$ are both of the same sign iff the left-hand side of (2.5) is positive, and this
happens iff $N(t_1)$ and $N(t_2)$ are of the same parity.
\medskip
Finally, one may consider the following problems. First, for how many $\g\;(\le T)$ the function $N(\g)$ is
even? Here the definition (1.5) is understood. Does one have
$$
\sum_{0<\g\le T, N(\g) = 2k}1 \;\sim\; \hf N(T)\qquad(T\to\infty)?
$$
This is certainly true if the conjecture that all zeros $\rho$ are simple is true. This, however,
seems to lie quite deep. The conjectures  {\it RH is true} and {\it all zeros are simple}
seem independent of one another. As far as it is known both conjectures could be true, or false, or one true
and the other one false.

\medskip
In the case when RH fails, it seems of interest to define, for a given $\g \,(>0)$,
$$
A(\g) := \sum_{\frac12<\b<1,\z(\b+i\g) = \z(\frac12+i\g)=0}1.
$$
It is clear that
$$
0 \,\le\, A(\g) \, \le\,  N(\g + \hf) - N(\g-\hf)  \,\ll\, \log\g.
$$
It is reasonable to expect that $A(\g) =0$ for almost all $\g$, but this is not easy to prove.

\vskip1cm
\medskip
\Refs
\medskip

\smallskip
\item{[1]} H.M. Edwards, Riemann's Zeta Function, Academic Press, New York,
1974.

\smallskip
\item{[2]} A. Ivi\'c, The Riemann zeta-function, John Wiley \&
Sons, New York, 1985 (reissue,  Dover, Mineola, New York, 2003).

\smallskip
\item{[3]} A. Ivi\'c,  Mean values of the Riemann zeta function, Tata Institute of Fund.
Research, LN's {\bf82}, (distributed by Springer Verlag, Berlin etc.), Bombay, 1991.
To be found online at
{\tt www.math.tifr.res.in/~publ/ln/tifr82.pdf}

\smallskip
\item{[4]} A. Ivi\'c,  On the multiplicity of  zeros of the zeta-function,
Bulletin CXVIII de l'Acad\'e\-mie Serbe des Sciences et des Arts - 1999,
Classe des Sciences math\'ematiques et naturelles,
Sciences math\'ematiques No. {\bf24}, pp. 119-131.

\smallskip
\item{[5]} A. Ivi\'c, The theory of Hardy's $Z$-function, Cambridge University Press,
Cambridge, 2013, 245 pp.

\smallskip
\item{[6]} A. Ivi\'c,  On the  multiplicites of zeros of $\z(s)$ and its values over short intervals,
Journal of Number Theory {\bf 185}(2018), 65-79.

\smallskip
\item{[7]} A.A. Karatsuba, Zero multiplicity and lower bound estimates of $|\z(s)|$,
Funct. Approx. Comment. Math. {\bf35}(2006), 195-207.

\smallskip

\item{[8]}
A. Selberg, On the remainder in the formula for $N(T)$, the number of zeros of $\z(s)$ in the
strip $0 < \g < T$, Avhandliger Norske  Videnkaps Akad. Oslo. I. Mat.-Naturv. Klasse, 1944, (1944).
no. 1, 27 pp, also in A. Selberg, Collected Works I, Springer, Berlin etc., 1989.

\smallskip
\item{[9]} E.C. Titchmarsh, The theory of the Riemann zeta-function, 2nd ed. edited by
D.R. Heath-Brown,  Clarendon Press, Oxford, 1986.

\smallskip
\item{[10]}
S. Wedeniwski, Results connected with the first 100 billion zeros of the Riemann
zeta function, 2002, at {\tt http://piologie.net/math/zeta.result.100billion.zeros.html}

\endRefs

\bigskip
\vskip1cm
\enddocument

\bye